\documentclass[11pt]{amsart}
\usepackage{amsmath}
\usepackage{amssymb,latexsym}
\usepackage{amscd}
\usepackage{xspace}
\usepackage{verbatim}
\begin{document}
\numberwithin{equation}{section}
\title[Equations with measures]{Remarks on nonlinear equations with measures}
\author{Moshe Marcus }
\address{Department of Mathematics, Technion\\
 Haifa 32000, ISRAEL}
 \email{marcusm@math.technion.ac.il}
\dedicatory{To the memory of  I. V. Skrypnik}
\date{\today}

\newcommand{\txt}[1]{\;\text{ #1 }\;}
\newcommand{\tbf}{\textbf}
\newcommand{\tit}{\textit}
\newcommand{\tsc}{\textsc}
\newcommand{\trm}{\textrm}
\newcommand{\mbf}{\mathbf}
\newcommand{\mrm}{\mathrm}
\newcommand{\bsym}{\boldsymbol}
\newcommand{\scs}{\scriptstyle}
\newcommand{\sss}{\scriptscriptstyle}
\newcommand{\txts}{\textstyle}
\newcommand{\dsps}{\displaystyle}
\newcommand{\fnz}{\footnotesize}
\newcommand{\scz}{\scriptsize}
\newcommand{\be}{\begin{equation}}
\newcommand{\bel}[1]{\begin{equation}\label{#1}}
\newcommand{\ee}{\end{equation}}
\newtheorem{subn}{\name}
\renewcommand{\thesubn}{}
\newcommand{\bsn}[1]{\def\name{#1$\!\!$}\begin{subn}}
\newcommand{\esn}{\end{subn}}
\newtheorem{sub}{\name}[section]
\newcommand{\dn}[1]{\def\name{#1}}   
\newcommand{\bs}{\begin{sub}}
\newcommand{\es}{\end{sub}}
\newcommand{\bsl}[1]{\begin{sub}\label{#1}}
\newcommand{\bth}[1]{\def\name{Theorem}\begin{sub}\label{t:#1}}
\newcommand{\blemma}[1]{\def\name{Lemma}\begin{sub}\label{l:#1}}
\newcommand{\bcor}[1]{\def\name{Corollary}\begin{sub}\label{c:#1}}
\newcommand{\bdef}[1]{\def\name{Definition}\begin{sub}\label{d:#1}}
\newcommand{\bprop}[1]{\def\name{Proposition}\begin{sub}\label{p:#1}}
\newcommand{\bnote}[1]{\def\name{\mdseries\scshape Notation}\begin{sub}\label{n:#1}}
\newcommand{\bproof}{\begin{proof}}
\newcommand{\eproof}{\end{proof}}
\newcommand{\bcom}{}
\newcommand{\req}{\eqref}
\newcommand{\rth}[1]{Theorem~\ref{t:#1}}
\newcommand{\rlemma}[1]{Lemma~\ref{l:#1}}
\newcommand{\rcor}[1]{Corollary~\ref{c:#1}}
\newcommand{\rdef}[1]{Definition~\ref{d:#1}}
\newcommand{\rprop}[1]{Proposition~\ref{p:#1}}
\newcommand{\rnote}[1]{Notation~\ref{n:#1}}
\newcommand{\BA}{\begin{array}}
\newcommand{\EA}{\end{array}}
\newcommand{\BAN}{\renewcommand{\arraystretch}{1.2}
\setlength{\arraycolsep}{2pt}\begin{array}}
\newcommand{\BAV}[2]{\renewcommand{\arraystretch}{#1}
\setlength{\arraycolsep}{#2}\begin{array}}
\newcommand{\BSA}{\begin{subarray}}
\newcommand{\ESA}{\end{subarray}}
\newcommand{\BAL}{\begin{aligned}}
\newcommand{\EAL}{\end{aligned}}
\newcommand{\BALG}{\begin{alignat}}
\newcommand{\EALG}{\end{alignat}}
\newcommand{\BALGN}{\begin{alignat*}}
\newcommand{\EALGN}{\end{alignat*}}
\newcommand{\note}[1]{\noindent\textit{#1.}\hspace{2mm}}
\newcommand{\Remark}{\note{Remark}}

\newcommand{\forevery}{\quad \forall}
\newcommand{\1}{\\[1mm]}
\newcommand{\2}{\\[2mm]}
\newcommand{\3}{\\[3mm]}
\newcommand{\set}[1]{\{#1\}}
\def\({{\rm (}}
\def\){{\rm )}}
\newcommand{\st}[1]{{\rm (#1)}}
\newcommand{\lra}{\longrightarrow}
\newcommand{\lla}{\longleftarrow}
\newcommand{\llra}{\longleftrightarrow}
\newcommand{\Lra}{\;\;\Longrightarrow\;\;}
\newcommand{\Lla}{\;\;\Longleftarrow\;\;}
\newcommand{\Llra}{\Longleftrightarrow}
\newcommand{\harpup}{\rightharpoonup}
\def\dar{\downarrow}
\def\uar{\uparrow}
\newcommand{\paran}[1]{\left (#1 \right )}
\newcommand{\sqrbr}[1]{\left [#1 \right ]}
\newcommand{\curlybr}[1]{\left \{#1 \right \}}
\newcommand{\abs}[1]{\left |#1\right |}
\newcommand{\norm}[1]{\left \|#1\right \|}
\newcommand{\angbr}[1]{\left< #1\right>}
\newcommand{\paranb}[1]{\big (#1 \big )}
\newcommand{\sqrbrb}[1]{\big [#1 \big ]}
\newcommand{\curlybrb}[1]{\big \{#1 \big \}}
\newcommand{\absb}[1]{\big |#1\big |}
\newcommand{\normb}[1]{\big \|#1\big \|}
\newcommand{\angbrb}[1]{\big\langle #1 \big \rangle}
\newcommand{\thkl}{\rule[-.5mm]{.3mm}{3mm}}
\newcommand{\thknorm}[1]{\thkl #1 \thkl\,}
\newcommand{\trinorm}[1]{|\!|\!| #1 |\!|\!|\,}
\newcommand{\vstrut}[1]{\rule{0mm}{#1}}
\newcommand{\rec}[1]{\frac{1}{#1}}
\newcommand{\opname}[1]{\mathrm{#1}\,}
\newcommand{\supp}{\opname{supp}}
\newcommand{\dist}{\opname{dist}}
\newcommand{\sign}{\opname{sign}}
\newcommand{\diam}{\opname{diam}}
\newcommand{\proj}{\opname{proj}}
\newcommand{\q}{\quad}
\newcommand{\qq}{\qquad}
\newcommand{\hsp}[1]{\hspace{#1mm}}
\newcommand{\vsp}[1]{\vspace{#1mm}}
\newcommand{\prt}{\partial}
\newcommand{\sms}{\setminus}
\newcommand{\ems}{\emptyset}
\newcommand{\ti}{\times}
\newcommand{\pr}{^\prime}
\newcommand{\ppr}{^{\prime\prime}}
\newcommand{\tl}{\tilde}
\newcommand{\wtl}{\widetilde}
\newcommand{\sbs}{\subset}
\newcommand{\sbeq}{\subseteq}
\newcommand{\nind}{\noindent}
\newcommand{\ovl}{\overline}
\newcommand{\unl}{\underline}
\newcommand{\nin}{\not\in}
\newcommand{\pfrac}[2]{\genfrac{(}{)}{}{}{#1}{#2}}
\newcommand{\tin}{\to\infty}
\newcommand{\ind}[1]{_{_{#1}}\!}
\newcommand{\chr}[1]{\chi\ind{#1}}
\newcommand{\rest}[1]{\big |\ind{#1}}
\newcommand{\num}[1]{{\rm (#1)}\hspace{2mm}}
\newcommand{\wkc}{weak convergence\xspace}
\newcommand{\wrto}{with respect to\xspace}
\newcommand{\cons}{consequence\xspace}
\newcommand{\consy}{consequently\xspace}
\newcommand{\Consy}{Consequently\xspace}
\newcommand{\Essy}{Essentially\xspace}
\newcommand{\essy}{essentially\xspace}
\newcommand{\mnz}{minimizer\xspace}
\newcommand{\sth}{such that\xspace}
\newcommand{\ngh}{neighborhood\xspace}
\newcommand{\nghs}{neighborhoods\xspace}
\newcommand{\seq}{sequence\xspace}
\newcommand{\sseq}{subsequence\xspace}
\newcommand{\locun}{locally uniformly\xspace}
\newcommand{\ifif}{if and only if\xspace}
\newcommand{\suff}{sufficiently\xspace}
\newcommand{\abc}{absolutely continuous\xspace}
\newcommand{\sol}{solution\xspace}
\newcommand{\subsol}{subsolution\xspace}
\newcommand{\supsol}{supersolution\xspace}
\newcommand{\Wlg}{Without loss of generality\xspace}
\newcommand{\wlg}{without loss of generality\xspace}
\newcommand{\bdw}{\partial\Gw}
\newcommand{\Capq}{C_{2/q,q'}}
\newcommand{\Cqq}{$C_{2/q,q'}$}
\newcommand{\finecl}{$C_{2/q,q'}$-finely closed \xspace}
\newcommand{\fineop}{$C_{2/q,q'}$-finely open \xspace}
\newcommand{\finetop}{$C_{2/q,q'}$-fine topology \xspace}
\newcommand{\Cqconv}{$C_{2/q,q'}$-convergent \xspace}

\newcommand{\Wq}{W^{2/q,q'}}
\newcommand{\Wqq}{W^{-2/q,q}}
\newcommand{\Wqdb}{W^{-2/q,q}_+(\bdw)}
\newcommand{\sbsq}{\overset{q}{\sbs}}
\newcommand{\smq}{\overset{q}{\sim}}
\newcommand{\app}[1]{\underset{#1}{\approx}}
\newcommand{\suppq}{\mathrm{supp}^q_{\bdw}\,}
\newcommand{\convq}{\overset{q}{\to}}
\newcommand{\barq}[1]{\bar{#1}^{^q}}
\newcommand{\prtq}{\partial_q}
\newcommand{\tr}{\mathrm{tr}}
\newcommand{\trV}{\mathrm{tr}\ind{V}}
\newcommand{\Tr}{\mathrm{Tr}\,}
\newcommand{\trR}{\mathrm{tr}\ind{\CR}}
\newcommand{\trin}[1]{\mathrm{tr}\ind{#1}}
\newcommand{\qcl}{$\Capq$-finely closed\xspace}
\newcommand{\qop}{$\Capq$-finely open\xspace}
\newcommand{\gsmod}{$\gs$-moderate\xspace}
\newcommand{\gsreg}{$\gs$-regular\xspace}
\newcommand{\qreg}{$q$-quasi regular\xspace}
\newcommand{\qeq}{$\Capq$-equivalent\xspace}
\newcommand{\ppf}{\underset{f}{\prec\prec}}
\newcommand{\ofrown}{\overset{\frown}}
\newcommand{\modcon}{\underset{mod}{\lra}}
\newcommand{\ugb}[1]{u\chr{\Gs_\gb(#1)}}
\newcommand{\mcon}{$q$-moderately convergent\xspace}
\newcommand{\mdiv}{$q$-moderately divergent\xspace}
\def\qsupp{q\text{-supp}\,}
\def\Lim{\,\text{\rm Lim}\,}
\def\muCR{\mu\ind{\CR}}
\def\vCR{v\ind{\CR}}
\def\bcom{}
\def\ga{\alpha}     \def\gb{\beta}       \def\gg{\gamma}
\def\gc{\chi}       \def\gd{\delta}      \def\ge{\epsilon}
\def\gth{\theta}                         \def\vge{\varepsilon}
\def\gf{\phi}       \def\vgf{\varphi}    \def\gh{\eta}
\def\gi{\iota}      \def\gk{\kappa}      \def\gl{\lambda}
\def\gm{\mu}        \def\gn{\nu}         \def\gp{\pi}
\def\vgp{\varpi}    \def\gr{\rho}        \def\vgr{\varrho}
\def\gs{\sigma}     \def\vgs{\varsigma}  \def\gt{\tau}
\def\gu{\upsilon}   \def\gv{\vartheta}   \def\gw{\omega}
\def\gx{\xi}        \def\gy{\psi}        \def\gz{\zeta}
\def\Gg{\Gamma}     \def\Gd{\Delta}      \def\Gf{\Phi}
\def\Gth{\Theta}
\def\Gl{\Lambda}    \def\Gs{\Sigma}      \def\Gp{\Pi}
\def\Gw{\Omega}     \def\Gx{\Xi}         \def\Gy{\Psi}

\def\CS{{\mathcal S}}   \def\CM{{\mathcal M}}   \def\CN{{\mathcal N}}
\def\CR{{\mathcal R}}   \def\CO{{\mathcal O}}   \def\CP{{\mathcal P}}
\def\CA{{\mathcal A}}   \def\CB{{\mathcal B}}   \def\CC{{\mathcal C}}
\def\CD{{\mathcal D}}   \def\CE{{\mathcal E}}   \def\CF{{\mathcal F}}
\def\CG{{\mathcal G}}   \def\CH{{\mathcal H}}   \def\CI{{\mathcal I}}
\def\CJ{{\mathcal J}}   \def\CK{{\mathcal K}}   \def\CL{{\mathcal L}}
\def\CT{{\mathcal T}}   \def\CU{{\mathcal U}}   \def\CV{{\mathcal V}}
\def\CZ{{\mathcal Z}}   \def\CX{{\mathcal X}}   \def\CY{{\mathcal Y}}
\def\CW{{\mathcal W}}
\def\BBA {\mathbb A}   \def\BBb {\mathbb B}    \def\BBC {\mathbb C}
\def\BBD {\mathbb D}   \def\BBE {\mathbb E}    \def\BBF {\mathbb F}
\def\BBG {\mathbb G}   \def\BBH {\mathbb H}    \def\BBI {\mathbb I}
\def\BBJ {\mathbb J}   \def\BBK {\mathbb K}    \def\BBL {\mathbb L}
\def\BBM {\mathbb M}   \def\BBN {\mathbb N}    \def\BBO {\mathbb O}
\def\BBP {\mathbb P}   \def\BBR {\mathbb R}    \def\BBS {\mathbb S}
\def\BBT {\mathbb T}   \def\BBU {\mathbb U}    \def\BBV {\mathbb V}
\def\BBW {\mathbb W}   \def\BBX {\mathbb X}    \def\BBY {\mathbb Y}
\def\BBZ {\mathbb Z}

\def\GTA {\mathfrak A}   \def\GTB {\mathfrak B}    \def\GTC {\mathfrak C}
\def\GTD {\mathfrak D}   \def\GTE {\mathfrak E}    \def\GTF {\mathfrak F}
\def\GTG {\mathfrak G}   \def\GTH {\mathfrak H}    \def\GTI {\mathfrak I}
\def\GTJ {\mathfrak J}   \def\GTK {\mathfrak K}    \def\GTL {\mathfrak L}
\def\GTM {\mathfrak M}   \def\GTN {\mathfrak N}    \def\GTO {\mathfrak O}
\def\GTP {\mathfrak P}   \def\GTR {\mathfrak R}    \def\GTS {\mathfrak S}
\def\GTT {\mathfrak T}   \def\GTU {\mathfrak U}    \def\GTV {\mathfrak V}
\def\GTW {\mathfrak W}   \def\GTX {\mathfrak X}    \def\GTY {\mathfrak Y}
\def\GTZ {\mathfrak Z}   \def\GTQ {\mathfrak Q}
\font\Sym= msam10
\def\SYM#1{\hbox{\Sym #1}}

\def\bmn{\mathbf{n}}
\def\bmm{\mathbf{m}}
\def\bma{\mathbf{a}}
\newcommand{\prtn}{\prt_{\mathbf{n}}}
\def\txin{\txt{in}}
\def\txon{\txt{on}}
\def\rhs{\text{right hand side \xspace}}
\def\lhs{\text{left hand side \xspace}}
\def\L1wsol{$L^1$-weak solution\xspace}
\def\W1p{W^{1,p}}
\def\Lr{L^1(\Gw,\gr)}
\def\lfs#1{\lfloor_{\sss  #1}}
\def\GTMr{\GTM(\Gw;\gr)}
\def\bgw{\bar\Gw}
\def\prtn{\prt_{\bmn}}
\def\loc{_{\mathrm{loc}}}
\def\bvp{boundary value problem\xspace}
\def\superh{superharmonic\xspace}
\def\subh{subharmonic\xspace}
\def\RN{\BBR^N}
\def\LVh{$L^V$ harmonic\xspace}
\def\LVsuph{$L^V$ superharmonic\xspace}
\def\LVmod{$L^V$ moderate}
\def\LVmodh{$L^V$ moderate harmonic}
\def\btr{boundary trace\xspace}
\def\Lip{Lipschitz\xspace}
\def\ssk{\smallskip}
\def\msk{\medskip}
\def\bsk{\bigskip}

\begin{abstract}
We study the Dirichlet boundary value  problem for equations with absorption of the form $-\Delta u+g\circ u=\mu$ in a bounded domain $\Omega\subset \BBR^N$ where $g$ is a continuous odd monotone increasing function.
Under some additional assumptions on $g$, we present necessary and sufficient conditions for existence when $\mu$ is a finite measure. We also discuss the notion of solution when the measure $\mu$ is positive and blows up on a compact subset of $\Gw$.
\end{abstract}
\maketitle
\section{Introduction}
In this paper we discuss some aspects of the  boundary value problem
\begin{equation}\label{gbvp}\BAL
  -\Gd u+g\circ u=\mu \q&\text{in }\Gw\\
  u=0\q&\text{on }\bdw,
  \EAL\end{equation}
where $\mu\in  \GTM_\gr(\Gw)$, i.e. $\mu$ is a Borel measure \sth
$$\int_\Gw\gr\,d|\mu|<\infty,  \q \gr(x)=dist(x,\bdw).$$
In addition we define a notion of solution in the case that $\mu$ is a positive Borel measure which may explode on a compact subset of the domain and discuss the question of existence and uniqueness in this case. We always assume that $g\in C(\BBR)$ is a monotone increasing function \sth $g(0)=0$. To simplify the presentation we also assume that $g$ is odd.

 A function $u\in L^1(\Gw)$ is a weak solution of the \bvp \req{gbvp},
$\mu\in \GTM_\gr$, if $u\in L^g_\gr(\Gw)$, i.e.
$$\int_\Gw\,g(u)\gr\,dx<\infty$$
and
\begin{equation}\label{inteq}
  \int_\Gw (-v\Gd\gf + g\circ v\,\gf)dx=\int_\Gw\gf\,d\mu
\end{equation}
for every $\gf\in C_0^2(\bgw)$ (= space of functions in $C^2(\bgw)$ vanishing on $\bdw$).

We say that $u$ is a solution of the equation
\begin{equation}\label{geq}
    -\Gd u+g\circ u=\mu \q \text{in } \Gw
\end{equation}
if $u$ and $g\circ u$ are in $L^1\loc(\Gw)$ and \req{inteq}  holds for every $\gf\in C_c^2(\Gw)$.

Brezis and Strauss \cite{BrS} proved that, if $\mu$ is an $L^1$ function the problem
possesses a unique solution.
This result does not extend to arbitrary measures in $\GTM_\gr(\Gw)$.

Denote by $\GTM^g_\gr$ the set of measures $\mu\in \GTM_\gr$ for which \req{gbvp} is solvable. A measure in $\GTM^g_\gr$ is called a \emph{$g$-good measure}.  It is known that, if a solution exists then it is unique.

 We say that $g$ is \emph{subcritical} if $\GTM^g_\gr=\GTM_\gr$.
Benilan and Brezis, \cite{Br70} and \cite{BeBr} proved that the following condition is sufficient for $g$ to
 be subcritical:
\begin{equation}\label{subcr}
\int_0^1 g(r^{2-N}))r^{N-1}dr<\infty.
\end{equation}
In the case that $g$ is a power non-linearity, i.e., $g=g_q$ where
$$g_q(t)=|t|^q\sign t  \q \text{in }\BBR,\q q>1,$$
this condition means that $q<q_c:=N/(N-2)$. Benilan and Brezis also proved that, if
$g=g_q$ and $q\geq q_c$,
 problem \req{gbvp} has no solution when $\mu$ is a Dirac measure.

 Later Baras and Pierre \cite{BP84} gave a complete characterization of $\GTM^g_\gr$ in the case that
$g=g_q$ with $q\geq q_c$. They proved that a finite measure $\mu$ is $g_q$-good if and only if $|\mu|$ does not charge
sets of $\bar C_{2,q'}$ capacity zero, $q'=q/(q-1)$. Here $\bar C_{\ga,p}$ denotes Bessel
capacity with the indicated indices.

In the present paper we extend the result of Baras and Pierre to a large class of non-linearities
and also discuss the notion of solution in the case that $\mu$ is a positive measure which  explodes on a
compact subset of $\Gw$.

\section{Statement of results}
Denote by $\CH$ the set  of even functions $h$ \sth
\begin{equation}\label{CG}\BAL
&h\in C^1(\BBR),\q h(0)=0, \q h\text{ is strictly convex,}\\
 & h'(0)=0,\q   h'(t)>0 \forevery t>0,\q\lim_{t\tin}h'(t)=\infty.
\EAL\end{equation}

For $h\in \CH$ denote by $L^h(\Gw)$ the corresponding Orlicz space in a domain $\Gw\sbs \RN$:
 $$L^h(\Gw)=\{f\in L^1\loc(\Gw)\mid\, \exists k>0: h\circ(f/k)\le 1\}$$
with the norm
$$\norm{f}\ind{L^h}=\inf\{k>0\mid\, h\circ(f/k)<\infty\}.$$
 Further denote by $h^*$ the conjugate of $h$. Since, by assumption, $h$ is strictly convex, $h'$ is strictly increasing so that,
$$h^*(t)=\int_0^t (h')^{-1}(s)ds.$$

Let $G$ be the Green kernel for $-\Gd$ in $\Gw$ and denote
$$\BBG_\mu(x)=\int_\Gw G(x,y)d\mu(y) \forevery x\in \Gw, \q \mu\in \GTM_\gr(\Gw).$$

For every $h\in\CH$, the capacity $C_{2,h}$ in $\Gw$ is defined as follows. For every compact set $E\sbs \Gw$ put:
\begin{equation}\label{C2g}
C_{2,h}(E)=\sup\{\mu(\Gw):\mu\in \GTM(\Gw),\;\mu\geq 0, \;\mu(E^c)=0,\; \norm{\BBG\mu}\ind{L^{h^*}}\le1\}.
\end{equation}
If $O$ is an open set:
$$C_{2,h}(O)=\sup\{C_{2,h}(E):\, E \sbs O,\; E\text{ compact.}\}$$
For an arbitrary set $A\sbs \Gw$ put
$$C_{2,h}(A)=\inf\{C_{2,h}(O):\, A\sbs O\sbs \Gw,\; O\text{ open.}\}$$
This definition is compatible with \req{C2g} : when $E$ is compact the value of $C_{2,h}(E)$ given by the above formula coincides with the value given by \req{C2g}, (see \cite{Or-cap}).

We say that $h$ satisfies the $\Gd_2$ condition if there exists $C>0$ \sth
$$h(a+b)\le c(h(a)+h(b)) \forevery a,b>0.$$
If $h\in \CH$ satisfies this condition then, $L^h$ is separable (see \cite{Kr-Rut}) and the capacity $C_{2,h}$ has the following additional properties (see \cite{Or-cap}).

Let $\Gw$ be a bounded domain in $\RN$. For every $A\sbs \Gw$,
\begin{equation}\label{C2g1}
   C_{2,h}(A)=\sup\{C_{2,h}(E):\, E \sbs A,\; E\text{ compact}\}
\end{equation}
and for every increasing \seq of sets $\{A_n\}$
 \begin{equation}\label{C2g2}
\lim C_{2,h}(A_n)=C_{2,h}(\cup\,A_n).
 \end{equation}
Furthermore, for every $A\sbs \Gw$
 \begin{equation}\label{C2g3}
  C_{2,h}(A)=\inf\{\norm{f}\ind{L^h}:\, f\in L^h(\Gw),\; \BBG_f\geq 1 \text{ on }A\}.
\end{equation}

If $h\in\CH$ and both $h$ and $h^*$ satisfy the $\Gd_2$ condition then $L^h$ is reflexive \cite{Kr-Rut}.

Finally we denote by $\CG$ the space of odd functions in $C(\BBR)$ \sth $h:=|g|\in \CH$ and by $\CG_2$ the set of functions $g\in \CG$ \sth $h$ and $h^*$ satisfy the $\Gd_2$ condition. For $g\in \CG$ put
$$L^g:=L^{|g|},\q C_{2,g}:=C_{2,|h|}, \q g^*(t)=|g|^*(t)\sign t \forevery t\in \BBR.$$


In the sequel we assume that $\Gw$ is a bounded domain of class $C^2$. The first theorem
provides a necessary and sufficient condition for the existence of a solution of \req{gbvp} in
 the spirit of \cite{BP84}.

\bth{Th.I}  Let $g\in \CG_2$ and let $\mu$ be a  measure in $ \GTM_\gr(\Gw)$.
Then problem \req{gbvp} possesses a solution
 if and only if $\mu$ vanishes on every compact set $E\sbs \Gw$ \sth $C_{2,g^*}(E)=0$. This condition will be indicated by the notation $\mu\prec C_{2,g^*}$.
\es

Next we consider problem \req{gbvp} when $ \mu$ is a positive Borel measure which may explode on a compact set $F\sbs \Gw$. In this part of the paper we assume that $g\in \CG_2$ and that $g$
satisfies the Keller -- Osserman condition \cite{Kell} and \cite{Oss}. This condition ensures that the set of solutions of
 \begin{equation}\label{geq0}
  -\Gd u+g\circ u=0
 \end{equation}
  in $\Gw$ is uniformly bounded in compact subsets of $\Gw$. Therefore, if $E\sbs \Gw$ and $E$  is   compact then there exists   a maximal solution of
\begin{equation}\label{max-sol}
 -\Gd u+g\circ u=0\q\text{in }\Gw\sms E, \q u=0\q\text{on }\bdw.
\end{equation}
This solution will be denoted by $U_E$.
\msk

\note{Notation} Consider the family of positive Borel measures $\mu$ in $\Gw$ \sth:

(1) There exists a compact set $F\sbs \Gw$ \sth, for every open set $O\supset F$,  $\mu(\Gw\sms\bar O)<\infty$ and

(2) $\mu(A)=\infty$ for every non-empty Borel set  $A\sbs F$.
\ssk

\noindent The set $F$ will be called the singular set of $\mu$. The family of measures $\mu$ of this type will be denoted by $\CB_\infty(\Gw)$.

\bdef{Definition 1}  Assume that $g\in \CG$ and that $g$ satisfies the Keller -- Osserman
condition. If $\nu\in \GTM^g_\gr(\Gw)$ denote by $v_\nu$ the solution of \req{gbvp} with $\mu$ replaced by $\nu$.

Let $\mu\in \CB_\infty(\Gw)$ and let $F$ be the singular set of $\mu$.
A function $u\in L^1\loc(\bgw\sms F)$ (i.e., $u\in L^1(\Gw\sms \bar O)$ for every \ngh $O$ of $F$) is a generalized solution  of \req{gbvp} if:

\num{i} $u$ satisfies \req{inteq}  for every $\gf\in C_0^2(\bgw)$
\sth $\supp\gf\sbs \Gw\sms F$.

\num{ii}$u\geq V_F:=\sup\{v_\nu:\,\nu\in \GTM_\gr^g(\Gw), \; \nu\geq 0, \; \supp\nu\sbs F\}.$
\es

\bth{Th.II} Assume that $g\in \CG_2$ and that $g$ satisfies the Keller -- Osserman condition.
Let $\mu\in \CB_\infty$ with singular set $F$.
Then:

\num{i} Problem \req{gbvp} has
a generalized solution if and only if $\mu$ vanishes on every compact set $E\sbs \Gw\sms F$ \sth
$C_{2,g^*}(E)=0$.

If $V_F=U_F$, where $V_F$ is defined as in \rdef{Definition 1} and $U_F$ is the maximal solution associated with $F$ (see \req{max-sol}) then the generalized solution is unique.


\num{ii} If $g$ satisfies the subcriticality condition \req{subcr} then problem \req{gbvp} possesses a unique generalized solution for every  $\mu\in \CB_\infty$.

\num{iii} Let $g=g_q$, $q\geq q_c$. If $\mu\prec C_{2,g^*}$ in $\Gw\sms F$ then \req{gbvp} possesses a unique solution.
\es

\section{Proof of \rth{Th.I}}
The proof is based on several lemmas. We assume throughout that the conditions of the theorem are satisfied.

Denote by $L^1_\gr(\Gw)$ the Lebesgue space  with weight $\gr$
 and by $L^g_\gr(\Gw)$ the Orlicz space with weight $\gr$.

Further denote by $W^kL^g(\Gw)$, $k\in \BBN$, the Orlicz-Sobolev space consisting of   functions $v\in L^g(\Gw)$ \sth $D^\ga v\in L^g(\Gw)$ for $|\ga|\le k$.

Under our assumptions the set of bounded functions in $L^g$ is dense in this space (see \cite{Kr-Rut}). \Consy, by \cite{DoTr},
$C^\infty(\bar \Gw)$ is dense in $W^kL^g(\Gw)$. As a consequence of the reflexivity of $L^g$ the  space $W^kL^g(\Gw)$ is reflexive. Let
$W^k_0L^g(\Gw)$ denote the closure of $C_c^\infty(\Gw)$ in $W^kL^g(\Gw)$. The dual of this space, denoted by $W^{-k}L^{g^*}(\Gw)$ is the linear hull of
$\{D^\ga f: f\in L^{g^*}(\Gw), \; |\ga\le k\}.$
The standard norm in $W^kL^g(\Gw)$ is given by
$$\norm{v}_{W^kL^g}=\sum_{|\ga|\le k} \norm{D^\ga v}\ind{L^g}$$
and the norm in $W^{-k}L^{g^*}$ is defined as the norm of the dual space of $W^{k}_0L^{g}$.

The spaces $W^kL^g_\gr$ and $W^{-k}L^{g^*}_\gr$ are defined in the same way.

\blemma{basics} If $\mu\in \GTM_\gr(\Gw)$ is a $g$-good measure then \req{gbvp} has a unique solution, which we denote by $v_\mu$.
The solution satisfies the inequality
\begin{equation}\label{vmu-est}
 \norm{v_\mu}\ind{L^1(\Gw)}+ \norm{v_\mu}\ind{L^g_\gr(\Gw)}\le C\norm{\mu}\ind{\GTM_\gr(\Gw)}
\end{equation}
where $C$ is a constant depending only on $g$ and $\Gw$.

If $\mu_j\in \GTM_\gr(\Gw)$, $j=1,2$ are $g$-good measures and $\mu_1\le \mu_2$ then $v_{\mu_1}\le v_{\mu_2}$.
\es

These results are well-known (see e.g. \cite{Vbook}).

\blemma{good1} Let $\mu\in \GTM_\gr(\Gw)$ be a positive measure \sth $\BBG_\mu\in L^g\loc(\Gw)$. Then $\mu$ is $g$ good.
\es
\proof Let $\{\Gw_n\}$ be a $C^2$ uniform exhaustion of $\Gw$. Then $\BBG_\mu\in L^g(\Gw_n)$ is a positive supersolution of problem \req{gbvp} in $\Gw_n$. Therefore -- as the zero function is a subsolution -- there exists a solution, say $u_n$, of \req{gbvp} in $\Gw_n$ and, by \rlemma{basics},
$$\int_{\Gw_n} u_ndx+\int_{\Gw_n} g\circ u_n\gr_n dx\leq C\int_{\Gw_n}\gr_n\,d\mu,$$
where $\gr_n(x)=\dist(x,\bdw_n)$ and $C$ is a constant depending only on $g$ and the $C^2$ character of $\Gw_n$. Since $\Gw_n\}$ is uniformly $C^2$, the constant may be chosen to be independent of $n$. Moreover  $\{u_n\}$ is increasing. Therefore $u=\lim u_n\in L^1(\Gw)\cap L^g_\gr(\Gw)$ is the solution of \req{gbvp}.
\qed

\blemma{admit} $(a)$  If $\mu\in \GTM_\gr$ and $|\mu|$ is $g$-good then $\mu$ is $g$-good.
$(b)$ $T\in W^{-2}L^g(\Gw)$ if and only if  $T=\Gd h$ for some $h\in L^g(\Gw)$.
$(c)$ If $\mu$ is a positive measure in  $W^{-2}L^g\loc(\Gw)$ then $\BBG_\mu\in L^g\loc(\Gw)$. If, in addition, $\mu\in \GTM_\gr(\Gw)$ then $\mu$ is $g$-good.
\es

\proof (a) Assuming that $|\mu|$ is $g$ -good, let $v$ be the solution of \req{gbvp} with $\mu$ replaced by $|\mu|$. Then $v$ is a supersolution and $-v$ is a subsolution of \req{gbvp}. Therefore \req{gbvp} has a solution.

(b)  If $T=\Gd h$ then, for every $\gf\in C_c^\infty(\Gw)$,
$$T(\gf)=\int_\Gw h\Gd \gf dx, \q |T(\gf)|\le \norm{h}\ind{L^g}\norm{\gf}\ind{W^{2}L^{g^*}}.$$
As $C_c^\infty$ is dense in $W^{2}_0L^{g^*}$, $T$ defines a continuous linear functional on this space; \consy $T\in W^{-2}L^g(\Gw)$.

On the other hand if $T\in W^{-2}L^g(\Gw)$, put
$$S(\Gd \gf):=T(\gf) \forevery \gf\in W^{2}_0L^{g^*}.$$
Note that for $\gf$ in this space we have $\gf=\BBG_{-\Gd\gf}$. Therefore $S$ is well defined on the subspace of $L^{g^*}$ given by
$\{\Gd\gf: \gf\in W^{2}_0L^{g^*}\}.$ Therefore there exists $h\in L^g(\Gw)$ \sth
$$T(\gf)=\int_\Gw h\Gd\gf\,dx \forevery \gf\in W^{2}_0L^{g^*}.$$
It follows that $T=\Gd h$.

(c) Let $\mu$ be a positive measure in $W^{-2}L^g\loc(\Gw)$. By part (b), if $\Gw'\Subset \Gw$ is a subdomain of class $C^2$ there exists  $h\in L^g(\Gw')$  \sth $\mu=\Gd h$. Then $h+\BBG_\mu$ is an  harmonic function in $\Gw'$;  \consy $\BBG_\mu\in L^g\loc(\Gw')$ and finally $\BBG_\mu\in L^g\loc(\Gw)$. If, in addition, $\mu\in \GTM_\gr(\Gw)$ then, by  \rlemma{good1}, $\mu$ is $g$ good.
\qed

\blemma{nesc1} Assume that $\mu\in \GTM_\gr(\Gw)$ is $g$ good. Then:

\num{i} There exists $f\in L^1_\gr(\Gw)$ and $\mu_0\in W^{-2}L^g\loc(\Gw)\cap\GTM_\gr(\Gw)$ \sth $\mu=f+\mu_0$.

\num{ii} $\mu\prec C_{2,g^*}$.
\es

\proof Assume that $\mu$ is $g$-good and let $u$ be the solution of \req{gbvp}. Then
$$\mu=f+\mu_0\txt{where} f:=g\circ u\in L^1_\gr,\;\mu_0:=\mu-g\circ u$$
and $u=\BBG_{\mu_0}\in L^g_\gr(\Gw)$. This implies that
$$\gf\mapsto \int_\Gw \gf\,d\mu_0=\int_\Gw \Gd\gf udx \forevery \gf\in C_c^\infty(\Gw)$$
is continuous on $C_0^2(\bar\Gw)$ \wrto the norm of $W^{2}L^{g^*}_\gr(\Gw)$. Therefore, the functional can be extended to a continuous linear functional on
$W^{2}L^{g^*}(\Gw')$ for every $\Gw'\Subset\Gw$. Thus $\mu_0\in W^{-2}L^g\loc(\Gw)\cap\GTM_\gr(\Gw)$.
\qed

(ii)  In view of \req{C2g1} it is sufficient to prove that $\mu$  vanishes on compact sets $E$ \sth $C_{2,g^*}(E)=0$.
\ssk

\note{Assertion} \emph{If $\nu\in W^{-2}L^{g}\loc(\Gw)$ then $\nu(E)=0$ for every compact set $E$ \sth  $C_{2,g^*}(E)=0$.}
\ssk

This assertion and part (i)  imply part (ii).

Suppose that there exists a  set $E$  \sth $C_{2,g^*}(E)=0$ and $\nu(E)\neq 0$. Then there exists a compact subset of $E$ on which $\nu$ has constant sign. Therefore we may assume that $E$ is compact and that $\nu$ is positive on $E$. We may assume that $\nu\in W^{-2}L^{g}(\Gw)$; otherwise we replace $\Gw$ by a $C^2$ domain $\Gw'\Subset\Gw$.

Let $\{V_n\}$ be a \seq of open \ngh{s} of $E$ \sth $\bar V_{n+1}\sbs V_n$ and $V_n\dar E$. Then there exists a \seq $\{\vgf_n\}$ in $C_c^\infty(\Gw)$ \sth $0\le \vgf_n\le 1$, $\vgf_n=1$ in $V_{n+1}$, $\supp\vgf_n\sbs V_n$ and  $\norm{\vgf_n}_{g^*}\to 0$.

This is proved in the same way as in the case of Bessel capacities. We use \req{C2g3} and the fact that $C^\infty(\bar\Gw)$ is dense in $W^2L^g_\gr(\Gw)$ \cite{DoTr}). Furthermore we use an extension of the  lemma on smooth truncation  \cite[Theorem 3.3.3]{AH}   to Sobolev-Orlicz spaces with an integral number of derivatives.  The extension is  straightforward.

Hence,
\begin{equation}\label{mu-g1}
\int_\Gw\vgf_n\,d\nu \to 0.
\end{equation}
On the other hand,
$$\int_\Gw \vgf_n d\nu\geq \nu(\bar V_{n+1})-|\nu|(V_n\sms \bar V_{n+1})\to \nu(E)>0.$$
This contradiction proves the assertion.
\qed

\msk

\blemma{mu<C} Let $\mu$ be a positive measure in $\GTM_\gr(\Gw)$. If $\mu$ vanishes on every compact set $E\sbs \Gw$ \sth $C_{2,g^*}(E)=0$ then $\mu$ is the limit of an increasing \seq of positive measures $\{\mu_n\}\sbs W^{-2}L^g(\Gw)$.
\es
\proof  Since $\mu$ is the limit of an increasing \seq of measures in $\GTM(\Gw)$ it is sufficient to prove the lemma for $\mu\in \GTM(\Gw)$.
Let $\vgf\in  W^{2}_0L^{g^*}(\Gw)$ and denote
$$\tl\vgf=\BBG_{\Gd \vgf}.$$
Then $\tl\vgf$ is equivalent to $\vgf$.

Suppose that $\{\vgf_n\}$ converges to $\vgf$ in $W^{2}_0L^{g^*}(\Gw)$. Then $\Gd\vgf_n\to\Gd\vgf$ in $L^{g^*}$. \Consy, by \cite[Theorem 4]{Or-cap},  there exists a  \sseq \sth $\tl\vgf_{n'}\to\tl\vgf$ $C_{2,g^*}$-a.e. (i.e., everywhere with the possible exception of a set of $C_{2,g^*}$-capacity zero). As $\mu$ vanishes on sets of capacity zero, it follows that $\tl\vgf_{n'}\to\tl\vgf$ $\mu$-a.e.{}.

Every $\vgf\in W^{2}_0L^{g^*}(\Gw)$  is the limit of a \seq $\{\vgf_n\}\sbs C^\infty_c(\Gw)$. Hence $\vgf_n\to \tl\vgf$ $\mu$-a.e. and \consy $\tl\vgf$ is $\mu$-measurable.

Therefore the functional $p:W^{2}_0L^{g^*}(\Gw)\mapsto [0,\infty]$ given by
$$p(\vgf):=\int_\Gw (\tl\vgf)_+d\mu$$
is well defined. The functional is sublinear, convex and l.s.c.: if $\vgf_n\to\vgf$  in $W^{2}_0L^{g^*}(\Gw)$ then (by Fatou's lemma)
$$p(\vgf)\le\liminf p(\vgf_n).$$
Furthermore,
$$p(a\vgf)=ap(\vgf) \forevery a>0.$$
Therefore the result follows by an application of the Hahn-Banach theorem, in the same way as in \cite[Lemma 4.2]{BP84}.
\qed
\msk

\note{Proof of \rth{Th.I}} By \rlemma{nesc1} the condition $\mu\prec C_{2,g^*}$ is necessary for the existence of a solution. We show that the condition is sufficient. 

If $\mu\prec C_{2,g^*}$ then $|\mu|\prec C_{2,g^*}$. By \rlemma{admit} if $|\mu|$ is $g$-good then $\mu$ is $g$-good. Therefore it remains  to prove the sufficiency of the condition for positive $\mu$. In this case, by \rlemma{mu<C}, there exists an increasing \seq of positive measures $\{\mu_n\}\sbs W^{-2}L^g(\Gw)$ \sth $\mu_n\uparrow\mu$. By \rlemma{admit} the measures $\mu_n$ are $g$-good. Denote by $u_n$ the solution of \req{gbvp} with $\mu$ replaced by $\mu_n$. By \rlemma{basics}, $u_n\geq 0$, $\{u_n\}$ increases and $\{u_n\}$ is bounded in $L^1(\Gw)\cap L^g_\gr(\Gw)$. Therefore  $u=\lim u_n\in L^1(\Gw)\cap L^g_\gr(\Gw)$ and $u_n\to u$ in this space. \Consy $u$ is the solution of \req{gbvp}.
\qed

\section{Proof of \rth{Th.II}}

(i) Let $\{O_n\}$ be a decreasing \seq of open sets \sth $\bar O_{n+1}\sbs O_n$, $\bar O_n\sbs \Gw$ and $O_n\dar F$ and $O_n$ is of class $C^2$. By \rth{Th.I}, the condition $\mu\prec C_{2,g^*}$ in $\Gw\sms F$ is necessary and sufficient for the existence of a solution of the equation
\begin{equation}\label{O_n}
  -\Gd u+g\circ u=\mu\q \text{in }\Gw_n:=\Gw\sms \bar O_n
\end{equation}
\sth $u=0$ on the boundary. By a standard argument, it follows that, under this condition:
for every $f\in L^1(\bdw\cup\prt O_n)$, \req{O_n} has a solution \sth $u=f$ on the boundary.
As $g$ satisfies the Keller -- Osserman condition, it also follows that \req{O_n} has a solution $u_n$ \sth $u_n=0$ on $\bdw$ and $u_n=\infty$ on $\prt O_n$.
Denote by $v_n$  the solution of \req{O_n} vanishing on $\bdw\cup\prt O_n$ and put
$$v_{0,\mu}=\lim v_n, \q \bar u_\mu=\lim u_n.$$
Then $v_{0,\mu}$ is the smallest positive solution of \req{O_n} vanishing on $\bdw$ while $\bar u_\mu$ is the largest such solution. In particular $\bar u_\mu\geq v_\nu$ for every $\nu\in \GTM^g_\gr$ \sth $\supp\nu\sbs F$. Thus $\bar u_\mu$ is the largest generalized solution of \req{gbvp}.

Next we construct the minimal  generalized solution of \req{gbvp}. The function $u_{0,\mu}+V_F$ is a supersolution and $\max(u_{0,\mu},V_F)$ is a subsolution of \req{O_n}, both vanishing on the boundary. Let $w_n$ denote the solution of \req{O_n} \sth $w_n=0$ on $\bdw$ and $w_n=\max(u_{0,\mu},V_F)$ on $\prt O_n$. Then
$$w_{n+1}\le w_n\le u_{0,\mu}+V_F$$
 and \consy, $w=\lim w_n$ is the smallest solution of \req{O_n} \sth
$$\max(u_{0,\mu},V_F)\le w\le u_{0,\mu}+V_F.$$
It follows that $w$ is a generalized solution of \req{gbvp}. Since any such solution dominates $\max(u_{0,\mu},V_F)$ it follows that $w$ is the smallest generalized solution of the problem. It is easy to see that $w=\unl{u}_\mu$ as given by  \req{max-sol}.

Since $g$ is convex, monotone increasing and $g(0)=0$ we have
$$g(a)+g(b)\le g(a+b) \forevery a,b\in \BBR_+.$$
Therefore $\bar u_\mu-u_{0,\mu}$ is a subsolution of \req{geq0} in $\Gw\sms F$. \Consy $\bar u_\mu-u_{0,\mu}\le U_F $ and
\begin{equation}\label{bar u}
 \max(u_{0,\mu},U_F)\le \bar u_\mu\le u_{0,\mu}+U_F.
\end{equation}

Put $\Gw_n=\Gw\sms \bar O_n$. Let $\unl{u}_n$ be the solution of the problem

$$\BAL
  -&\Gd u+g\circ u=\mu\q \text{in }\Gw_n, \\ &u=V_F \;\text{ on }\prt O_n,\q u=0\;\text{ on } \bdw.
  \EAL$$
Then $\{\unl u_n\}$ increases and $\unl u=\lim \unl u_n$.

Similarly, if $\bar u_n$ is the solution of the problem

$$\BAL
  -&\Gd u+g\circ u=\mu\q \text{in }\Gw_n, \\ &u=U_F \;\text{ on }\prt O_n,\q u=0\;\text{ on } \bdw.
  \EAL$$
then $\{\bar u_n\}$ increases and, in view of \req{bar u}, $\bar u=\lim \bar u_n$.
Therefore, if $V_F=U_F$ then $\bar u_\mu=\unl{u}_\mu$.

(ii) We assume that in addition to the other conditions of the theorem,  $g$ satisfies the subcriticality condition. In this case, for every point $z\in \Gw$ and $k\in \BBR$,  there exists a  solution $u_{k,z}$ of the problem
\begin{equation}\label{gdz}
  -\Gd u+g\circ u=k\gd_z \q\text{in }\Gw,\q u=0 \q \text{on }\bdw.
\end{equation}

Put $w_z=\lim_{k\tin}u_{k,z}$. By definition $w_z=V_{\{z\}}$. We also have $w_z=U_{\{z\}}$. This follows from the fact that $g$ satisfies the Keller -- Osserman condition. This condition implies that there exists a decreasing function $\psi\in C(0,\infty)$ \sth  $\psi(t)\tin$ as $t\to 0$ and every positive solution $u$ of \req{gdz} satisfies
$$C_2\psi(|x-z|) \le u(x) \le C_1\psi(|x-z|)$$.The constant $C_1$ depends only on $g,N$. Because of the boundary condition the constant  $C_2$ depends on $z$. However for $z$ in a compact subset of $\Gw$ one can choose $C_2$ to be independent of $z$.

This inequality implies that
$$w_z\le U_{\{z\}}\le C_1/C_2 w_z.$$
If $F$ is a compact subset of $\Gw$ put
$$F'=\{x\in \Gw: dist(x,F)\le \rec{2}\dist(F,\bdw)\}.$$

Let $x\in F'\sms F$ and let $z$ be a point in $F$ \sth $|x-z|=\dist(x,F)$. Then there exists a positive constant $C(F)$ \sth
$$ C(F)\psi(|x-z|)\le U_z(x)\le V_F(x)\le U_F(x)\le C_1\psi(|x-z|).$$
It follows that there exists a constant $c$ \sth
\begin{equation}\label{U<cV}
 U_F(x)\le cV_F(x)
\end{equation}
for every $x\in F'$. Since $U_F$ and $V_F$ vanish on $\bdw$ it follows that \req{U<cV} (with possibly a larger constant) remains valid in $\Gw\sms F'$. This is verified by a standard argument using Harnack's inequality and the fact that $g$ satisfies the Keller -- Osserman condition. Thus \req{U<cV} is valid in $\Gw\sms F$. By an argument similar to the one introduced in \cite[Theorem 5.4]{MVsub}, this inequality implies that $U_F=V_F$.

(iii) For the case considered  here, it was proved in  \cite{MVcapest}  that $U_F=V_F$. Therefore uniqueness follows from part (i).
\qed

\end{document}